\documentclass{amsart}
\usepackage{amsfonts}

\newtheorem{theorem}{Theorem}
\theoremstyle{plain}

\newtheorem{corollary}{Corollary}

\newtheorem{remark}{Remark}

\numberwithin{equation}{section}
\input{tcilatex}

\begin{document}
\title[Inequalities for the Norm and Numerical Radius]{Inequalities for the
Norm and Numerical Radius of Composite Operators in Hilbert Spaces}
\author{S.S. Dragomir}
\address{School of Computer Science and Mathematics\\
Victoria University of Technology\\
PO Box 14428, Melbourne City\\
Victoria 8001, Australia.}
\email{sever.dragomir@vu.edu.au}
\urladdr{http://rgmia.vu.edu.au/dragomir}
\date{5 September, 2005}

\begin{abstract}
Some new inequalities for the norm and the numerical radius of composite
operators generated by a pair of operators are given.
\end{abstract}

\maketitle

\section{Introduction}

Let $\left( H;\left\langle \cdot ,\cdot \right\rangle \right) $ be a complex
Hilbert space. The \textit{numerical range} of an operator $T$ is the subset
of the complex numbers $\mathbb{C}$ given by \cite[p. 1]{GR}:%
\begin{equation}
W\left( T\right) =\left\{ \left\langle Tx,x\right\rangle ,\ x\in H,\
\left\Vert x\right\Vert =1\right\} .  \label{1.1}
\end{equation}%
It is well known that (see \cite{GR}):

\begin{enumerate}
\item[(i)] The numerical range of an operator is convex;

\item[(ii)] The spectrum of an operator is contained in the closure of its
numerical range;

\item[(iii)] $T$ is self-adjoint if and only if $W\left( T\right) $ is real.
\end{enumerate}

The \textit{numerical radius} $w\left( T\right) $ of an operator $T$ on $H$
is defined by \cite[p. 8]{GR}%
\begin{equation}
w\left( T\right) :=\sup \left\{ \left\vert \lambda \right\vert ,\lambda \in
W\left( T\right) \right\} =\sup \left\{ \left\vert \left\langle
Tx,x\right\rangle \right\vert ,\left\Vert x\right\Vert =1\right\} .
\label{1.2}
\end{equation}

It is well known that $w\left( \cdot \right) $ is a norm on the Banach
algebra $B\left( H\right) $ of all bounded linear operators acting on $H$
and the following inequality holds true:%
\begin{equation}
w\left( T\right) \leq \left\Vert T\right\Vert \leq 2w\left( T\right) .
\label{1.3}
\end{equation}

We recall some classical results involving the numerical radius of two
linear operators $A,B.$

The following general result for the product of two operators holds \cite[p.
37]{GR}:

\begin{theorem}
\label{t1.1}If $A,B$ are two bounded linear operators on the Hilbert space $%
\left( H,\left\langle \cdot ,\cdot \right\rangle \right) ,$ then%
\begin{equation}
w\left( AB\right) \leq 4w\left( A\right) w\left( B\right) .  \label{1.4}
\end{equation}%
In the case that $AB=BA,$ then%
\begin{equation}
w\left( AB\right) \leq 2w\left( A\right) w\left( B\right) .  \label{1.5}
\end{equation}
\end{theorem}

The following results are also well known \cite[p. 38]{GR}.

\begin{theorem}
\label{t1.2}If $A$ is a unitary operator that commutes with another operator 
$B,$ then%
\begin{equation}
w\left( AB\right) \leq w\left( B\right) .  \label{1.6}
\end{equation}%
If $A$ is an isometry and $AB=BA,$ then (\ref{1.6}) also holds true.
\end{theorem}

We say that $A$ and $B$ \textit{double commute} if $AB=BA$ and $AB^{\ast
}=B^{\ast }A.$

The following result holds \cite[p. 38]{GR}.

\begin{theorem}[Double commute]
\label{t1.3}If the operators $A$ and $B$ double commute, then%
\begin{equation}
w\left( AB\right) \leq w\left( B\right) \left\Vert A\right\Vert .
\label{1.7}
\end{equation}
\end{theorem}

As a consequence of the above, we have \cite[p. 39]{GR}:

\begin{corollary}
\label{c1.1}Let $A$ be a normal operator commuting with $B.$ Then%
\begin{equation}
w\left( AB\right) \leq w\left( A\right) w\left( B\right) .  \label{1.8}
\end{equation}
\end{corollary}

For other results and historical comments on the above see \cite[p. 39--41]%
{GR}. For more results on the numerical radius, see \cite{H}.

The main aim of this paper is to establish some new inequalities for
composite operators generated by a pair of operators $\left( A,B\right) $
under various assumptions. Namely, in one side, several inequalities
involving the norm 
\begin{equation*}
\left\Vert \frac{A^{\ast }A+B^{\ast }B}{2}\right\Vert
\end{equation*}%
and the numerical radius $w\left( B^{\ast }A\right) $ are established. On
the other side, upper bounds for the nonnegative quantities 
\begin{equation*}
\left\Vert A\right\Vert \left\Vert B\right\Vert -w\left( B^{\ast }A\right) 
\text{ and }\left\Vert A\right\Vert ^{2}\left\Vert B\right\Vert
^{2}-w^{2}\left( B^{\ast }A\right)
\end{equation*}%
under special conditions for the operators involved are also given.

Applications for normal operators are provided as well.

\section{The Results}

The following result may be stated:

\begin{theorem}
\label{t2.1}Let $A,B:H\rightarrow H$ be two bounded linear operators on the
Hilbert space $\left( H,\left\langle \cdot ,\cdot \right\rangle \right) .$
If $r>0$ and%
\begin{equation}
\left\Vert A-B\right\Vert \leq r,  \label{2.1}
\end{equation}%
then%
\begin{equation}
\left\Vert \frac{A^{\ast }A+B^{\ast }B}{2}\right\Vert \leq w\left( B^{\ast
}A\right) +\frac{1}{2}r^{2}.  \label{2.2}
\end{equation}
\end{theorem}

\begin{proof}
For any $x\in H,$ $\left\Vert x\right\Vert =1,$ we have from (\ref{2.1}) that%
\begin{equation}
\left\Vert Ax\right\Vert ^{2}+\left\Vert Bx\right\Vert ^{2}\leq 2\func{Re}%
\left\langle Ax,Bx\right\rangle +r^{2}.  \label{2.3}
\end{equation}%
However%
\begin{align*}
\left\Vert Ax\right\Vert ^{2}+\left\Vert Bx\right\Vert ^{2}& =\left\langle
\left( A^{\ast }A\right) x,x\right\rangle +\left\langle \left( B^{\ast
}B\right) x,x\right\rangle \\
& =\left\langle \left( A^{\ast }A+B^{\ast }B\right) x,x\right\rangle
\end{align*}%
and by (\ref{2.3}) we obtain%
\begin{equation}
\left\langle \left( A^{\ast }A+B^{\ast }B\right) x,x\right\rangle \leq
2\left\vert \left\langle \left( B^{\ast }A\right) x,x\right\rangle
\right\vert +r^{2}  \label{2.4}
\end{equation}%
for any $x\in H,$ $\left\Vert x\right\Vert =1.$

Taking the supremum over $x\in H,$ $\left\Vert x\right\Vert =1$ in (\ref{2.4}%
) we get%
\begin{equation}
w\left( A^{\ast }A+B^{\ast }B\right) \leq 2w\left( B^{\ast }A\right) +r^{2}
\label{2.5}
\end{equation}%
and since the operator $A^{\ast }A+B^{\ast }B$ is self-adjoint, hence 
\begin{equation*}
w\left( A^{\ast }A+B^{\ast }B\right) =\left\Vert A^{\ast }A+B^{\ast
}B\right\Vert
\end{equation*}%
and by (\ref{2.5}) we deduce the desired inequality (\ref{2.2}).
\end{proof}

\begin{remark}
We observe that, from the proof of the above theorem, we have the
inequalities%
\begin{equation}
0\leq \left\Vert \frac{A^{\ast }A+B^{\ast }B}{2}\right\Vert -w\left( B^{\ast
}A\right) \leq \frac{1}{2}\left\Vert A-B\right\Vert ^{2},  \label{2.6}
\end{equation}%
provided that $A,B$ are bounded linear operators in $H.$

The second inequality in (\ref{2.6}) is obvious while the first inequality
follows by the fact that%
\begin{eqnarray*}
\left\langle \left( A^{\ast }A+B^{\ast }B\right) x,x\right\rangle
&=&\left\Vert Ax\right\Vert ^{2}+\left\Vert Bx\right\Vert ^{2} \\
&\geq &2\left\Vert Ax\right\Vert \left\Vert Bx\right\Vert \geq 2\left\vert
\left\langle \left( B^{\ast }A\right) x,x\right\rangle \right\vert
\end{eqnarray*}%
for any $x\in H.$
\end{remark}

The inequality (\ref{2.2}) is obviously a reach source of particular
inequalities of interest.

Indeed, if we assume, for $\lambda \in \mathbb{C}$ and a bounded linear
operator $T,$ that we have 
\begin{equation}
\left\Vert T-\lambda T^{\ast }\right\Vert \leq r,  \label{2.6.1}
\end{equation}%
for a given positive number $r,$ then by (\ref{2.6}) we deduce the inequality%
\begin{equation}
0\leq \left\Vert \frac{T^{\ast }T+\left\vert \lambda \right\vert
^{2}TT^{\ast }}{2}\right\Vert -\left\vert \lambda \right\vert w\left(
T^{2}\right) \leq \frac{1}{2}r^{2}.  \label{2.6.2}
\end{equation}%
Moreover, if $T$ is assumed to be \textit{normal}, i.e., $T^{\ast
}T=TT^{\ast },$ then by (\ref{2.6.2}) we recapture our result from \cite%
{SSD1}, namely:%
\begin{equation*}
0\leq \frac{1+\left\vert \lambda \right\vert ^{2}}{2}\left\Vert T\right\Vert
^{2}-\left\vert \lambda \right\vert w\left( T^{2}\right) \leq \frac{1}{2}%
r^{2}
\end{equation*}%
provided the normal operator $T$ satisfies (\ref{2.6.1}).

Now, if we assume that for $\lambda \in \mathbb{C}$ and a bounded linear
operator $V$ we have that%
\begin{equation}
\left\Vert V-\lambda I\right\Vert \leq r,  \label{2.6.3}
\end{equation}%
where $I$ is the identity operator on $H,$ then by (\ref{2.2}) we deduce the
inequality%
\begin{equation*}
0\leq \left\Vert \frac{V^{\ast }V+\left\vert \lambda \right\vert ^{2}I}{2}%
\right\Vert -\left\vert \lambda \right\vert w\left( V\right) \leq \frac{1}{2}%
r^{2}.
\end{equation*}

As a dual approach, the following result may be noted as well:

\begin{theorem}
\label{t2.2}Let $A,B:H\rightarrow H$ be two bounded linear operators on the
Hilbert space $H.$ Then%
\begin{equation}
\left\Vert \frac{A+B}{2}\right\Vert ^{2}\leq \frac{1}{2}\left[ \left\Vert 
\frac{A^{\ast }A+B^{\ast }B}{2}\right\Vert +w\left( B^{\ast }A\right) \right]
.  \label{2.7}
\end{equation}
\end{theorem}

\begin{proof}
We obviously have%
\begin{align*}
\left\Vert Ax+Bx\right\Vert ^{2}& =\left\Vert Ax\right\Vert ^{2}+2\func{Re}%
\left\langle Ax,Bx\right\rangle +\left\Vert Bx\right\Vert ^{2} \\
& \leq \left\langle \left( A^{\ast }A+B^{\ast }B\right) x,x\right\rangle
+2\left\vert \left\langle \left( B^{\ast }A\right) x,x\right\rangle
\right\vert
\end{align*}%
for any $x\in H.$

Taking the supremum over $x\in H,$ $\left\Vert x\right\Vert =1,$ we get%
\begin{align*}
\left\Vert A+B\right\Vert ^{2}& \leq w\left( A^{\ast }A+B^{\ast }B\right)
+2w\left( B^{\ast }A\right) \\
& =\left\Vert A^{\ast }A+B^{\ast }B\right\Vert +2w\left( B^{\ast }A\right) ,
\end{align*}%
from where we get the desired inequality (\ref{2.7}).
\end{proof}

\begin{remark}
The inequality (\ref{2.7}) can generate some interesting particular results
such as the following inequality%
\begin{equation}
\left\Vert \frac{T+T^{\ast }}{2}\right\Vert ^{2}\leq \frac{1}{2}\left[
\left\Vert \frac{T^{\ast }T+TT^{\ast }}{2}\right\Vert +w\left( T^{2}\right) %
\right] ,  \label{2.7.1}
\end{equation}%
holding for each bounded linear operator $T:H\rightarrow H.$

If, in particular, $T$ is assumed to be normal, then (\ref{2.7.1}) becomes%
\begin{equation}
\left\Vert \frac{T+T^{\ast }}{2}\right\Vert ^{2}\leq \frac{1}{2}\left[
\left\Vert T\right\Vert ^{2}+w\left( T^{2}\right) \right] \left( \leq
\left\Vert T\right\Vert ^{2}\right) .  \label{2.7.2}
\end{equation}%
Its is well known that for $V$ a bounded linear operator on $H$, by the
convexity property of the $\left\Vert \cdot \right\Vert ^{2}$ on $B\left(
H\right) $ we have%
\begin{equation}
\left\Vert \frac{V+V^{\ast }}{2}\right\Vert ^{2}\leq \frac{\left\Vert
V\right\Vert ^{2}+\left\Vert V^{\ast }\right\Vert ^{2}}{2}=\left\Vert
V\right\Vert ^{2}.  \label{2.7.3}
\end{equation}%
Therefore, the inequality (\ref{2.7.2}), which holds for normal operators,
is a refinement of (\ref{2.7.3}) that holds for any bounded linear operator.
This result has been obtained in a different manner in the earlier paper 
\cite{SSD1}.
\end{remark}

The following result may be stated as well.

\begin{theorem}
\label{t2.3}Let $A,B:H\rightarrow H$ be two bounded linear operators on the
Hilbert space $H$ and $p\geq 2.$ Then%
\begin{equation}
\left\Vert \frac{A^{\ast }A+B^{\ast }B}{2}\right\Vert ^{\frac{p}{2}}\leq 
\frac{1}{4}\left[ \left\Vert A-B\right\Vert ^{p}+\left\Vert A+B\right\Vert
^{p}\right] .  \label{2.8}
\end{equation}
\end{theorem}

\begin{proof}
We use the following inequality for vectors in inner product spaces obtained
by Dragomir and S\'{a}ndor in \cite{DS}:%
\begin{equation}
2\left( \left\Vert a\right\Vert ^{p}+\left\Vert b\right\Vert ^{p}\right)
\leq \left\Vert a+b\right\Vert ^{p}+\left\Vert a-b\right\Vert ^{p}
\label{2.9}
\end{equation}%
for any $a,b\in H$ and $p\geq 2.$

Utilising (\ref{2.9}) we may write%
\begin{equation}
2\left( \left\Vert Ax\right\Vert ^{p}+\left\Vert Bx\right\Vert ^{p}\right)
\leq \left\Vert Ax+Bx\right\Vert ^{p}+\left\Vert Ax-Bx\right\Vert ^{p}
\label{2.10}
\end{equation}%
for any $x\in H.$

Now, observe that%
\begin{equation*}
\left\Vert Ax\right\Vert ^{p}+\left\Vert Bx\right\Vert ^{p}=\left(
\left\Vert Ax\right\Vert ^{2}\right) ^{\frac{p}{2}}+\left( \left\Vert
Bx\right\Vert ^{2}\right) ^{\frac{p}{2}}
\end{equation*}%
and by the elementary inequality:%
\begin{equation*}
\frac{\alpha ^{q}+\beta ^{q}}{2}\geq \left( \frac{\alpha +\beta }{2}\right)
^{q},\quad \alpha ,\beta \geq 0\text{ \ and \ }q\geq 1
\end{equation*}%
we have%
\begin{align}
\left( \left\Vert Ax\right\Vert ^{2}\right) ^{\frac{p}{2}}+\left( \left\Vert
Bx\right\Vert ^{2}\right) ^{\frac{p}{2}}& \geq 2^{1-\frac{p}{2}}\left(
\left\Vert Ax\right\Vert ^{2}+\left\Vert Bx\right\Vert ^{2}\right) ^{\frac{p%
}{2}}  \label{2.11} \\
& =2^{1-\frac{p}{2}}\left[ \left\langle \left( A^{\ast }A+B^{\ast }B\right)
x,x\right\rangle \right] ^{\frac{p}{2}}.  \notag
\end{align}%
Combining (\ref{2.10}) with (\ref{2.11}) we get%
\begin{equation}
\frac{1}{4}\left[ \left\Vert Ax-Bx\right\Vert ^{p}+\left\Vert
Ax+Bx\right\Vert ^{p}\right] \geq \left\vert \left\langle \left( \frac{%
A^{\ast }A+B^{\ast }B}{2}\right) x,x\right\rangle \right\vert ^{\frac{p}{2}}
\label{2.12}
\end{equation}%
for any $x\in H,$ $\left\Vert x\right\Vert =1.$ Taking the supremum over $%
x\in H,$ $\left\Vert x\right\Vert =1,$ and taking into account that 
\begin{equation*}
w\left( \frac{A^{\ast }A+B^{\ast }B}{2}\right) =\left\Vert \frac{A^{\ast
}A+B^{\ast }B}{2}\right\Vert ,
\end{equation*}%
we deduce the desired result (\ref{2.8}).
\end{proof}

\begin{remark}
If $p=2,$ then we have the inequality:%
\begin{equation}
\left\Vert \frac{A^{\ast }A+B^{\ast }B}{2}\right\Vert \leq \left\Vert \frac{%
A-B}{2}\right\Vert ^{2}+\left\Vert \frac{A+B}{2}\right\Vert ^{2},
\label{2.13}
\end{equation}%
for any $A,B$ bounded linear operators. This result can also be obtained
directly on utilising the parallelogram identity.

We also should observe that for $A=T$ and $B=T^{\ast },$ $T$ a normal
operator, the inequality (\ref{2.8}) becomes%
\begin{equation*}
\left\Vert T\right\Vert ^{p}\leq \frac{1}{4}\left[ \left\Vert T-T^{\ast
}\right\Vert ^{p}+\left\Vert T+T^{\ast }\right\Vert ^{p}\right] ,
\end{equation*}%
where $p\geq 2.$
\end{remark}

The following result may be stated as well.

\begin{theorem}
\label{t2.4}Let $A,B:H\rightarrow H$ be two bounded linear operators on the
Hilbert space $H$ and $r\geq 1.$ If $A^{\ast }A\geq B^{\ast }B$ in the
operator order or, equivalently, $\left\Vert Ax\right\Vert \geq \left\Vert
Bx\right\Vert $ for any $x\in H,$ then:%
\begin{equation}
\left\Vert \frac{A^{\ast }A+B^{\ast }B}{2}\right\Vert ^{r}\leq \left\Vert
A\right\Vert ^{r-1}\left\Vert B\right\Vert ^{r-1}w\left( B^{\ast }A\right) +%
\frac{1}{2}r^{2}\left\Vert A\right\Vert ^{2r-2}\left\Vert A-B\right\Vert
^{2}.  \label{2.14}
\end{equation}
\end{theorem}

\begin{proof}
We use the following inequality for vectors in inner product spaces due to
Goldstein, Ryff and Clarke \cite{GRC}:%
\begin{equation}
\left\Vert a\right\Vert ^{2r}+\left\Vert b\right\Vert ^{2r}\leq 2\left\Vert
a\right\Vert ^{r-1}\left\Vert b\right\Vert ^{r-1}\func{Re}\left\langle
a,b\right\rangle +r^{2}\left\Vert a\right\Vert ^{2r-2}\left\Vert
a-b\right\Vert ^{2},  \label{2.15}
\end{equation}%
where $r\geq 1,$ $a,b\in H$ and $\left\Vert a\right\Vert \geq \left\Vert
b\right\Vert .$

Utilising (\ref{2.15}) we can state that:%
\begin{multline}
\quad \left\Vert Ax\right\Vert ^{2r}+\left\Vert Bx\right\Vert ^{2r}
\label{2.16} \\
\leq 2\left\Vert Ax\right\Vert ^{r-1}\left\Vert Bx\right\Vert
^{r-1}\left\vert \left\langle Ax,Bx\right\rangle \right\vert
+r^{2}\left\Vert Ax\right\Vert ^{2r-2}\left\Vert Ax-Bx\right\Vert ^{2},\quad
\end{multline}%
for any $x\in H.$

As in the proof of Theorem \ref{t2.3}, we also have%
\begin{equation}
2^{1-r}\left[ \left\langle \left( A^{\ast }A+B^{\ast }B\right)
x,x\right\rangle \right] ^{r}\leq \left\Vert Ax\right\Vert ^{2r}+\left\Vert
Bx\right\Vert ^{2r},  \label{2.17}
\end{equation}%
for any $x\in H.$

Therefore, by (\ref{2.16}) and (\ref{2.17}) we deduce%
\begin{multline}
\left[ \left\langle \left( \frac{A^{\ast }A+B^{\ast }B}{2}\right)
x,x\right\rangle \right] ^{r}  \label{2.18} \\
\leq \left\Vert Ax\right\Vert ^{r-1}\left\Vert Bx\right\Vert
^{r-1}\left\vert \left\langle Ax,Bx\right\rangle \right\vert +\frac{1}{2}%
r^{2}\left\Vert A\right\Vert ^{2r-2}\left\Vert Ax-Bx\right\Vert ^{2}
\end{multline}%
for any $x\in H.$

Taking the supremum in (\ref{2.18}) we obtain the desired result (\ref{2.14}%
).
\end{proof}

\begin{remark}
Following \cite[p. 156]{GR}, we recall that the bounded linear operator $V$
is hyponormal, if 
\begin{equation*}
\left\Vert V^{\ast }x\right\Vert \leq \left\Vert Vx\right\Vert \text{ for
all }x\in H.
\end{equation*}%
Now, if we choose in (\ref{2.14}) $A=V$ and $B=V^{\ast },$ then we get the
inequality%
\begin{equation}
\left\Vert \frac{V^{\ast }V+VV^{\ast }}{2}\right\Vert ^{r}\leq \left\Vert
V\right\Vert ^{2r-2}\left[ w\left( V^{2}\right) +\frac{1}{2}r^{2}\left\Vert
V-V^{\ast }\right\Vert ^{2}\right] ,  \label{2.18.1}
\end{equation}%
holding for any hyponormal operator $V$ and any $r\geq 1.$ In particular, if 
$V=T,$ a normal operator, then from (\ref{2.18.1}) we deduce for $r=1$ the
inequality%
\begin{equation}
\left\Vert T\right\Vert ^{2}\leq w\left( T^{2}\right) +\frac{1}{2}\left\Vert
T-T^{\ast }\right\Vert ^{2},  \label{2.18.2}
\end{equation}%
that has been pointed out in the earlier paper \cite{SSD1}.
\end{remark}

\section{Further Inequalities for an Invertible Operator}

In this section we assume that $B:H\rightarrow H$ is an invertible bounded
linear operator and let $B^{-1}:H\rightarrow H$ be its inverse. Then,
obviously,%
\begin{equation}
\left\Vert Bx\right\Vert \geq \frac{1}{\left\Vert B^{-1}\right\Vert }%
\left\Vert x\right\Vert \quad \text{for any \ }x\in H,  \label{3.1}
\end{equation}%
where $\left\Vert B^{-1}\right\Vert $ denotes the norm of the inverse $%
B^{-1}.$

The following result holds true:

\begin{theorem}
\label{t3.1}Let $A,B:H\rightarrow H$ be two bounded linear operators on $H$
and $B$ is invertible such that, for a given $r>0,$%
\begin{equation}
\left\Vert A-B\right\Vert \leq r.  \label{3.2}
\end{equation}%
Then:%
\begin{equation}
\left\Vert A\right\Vert \leq \left\Vert B^{-1}\right\Vert \left[ w\left(
B^{\ast }A\right) +\frac{1}{2}r^{2}\right] .  \label{3.3}
\end{equation}
\end{theorem}

\begin{proof}
The condition (\ref{3.2}) is obviously equivalent to:%
\begin{equation}
\left\Vert Ax\right\Vert ^{2}+\left\Vert Bx\right\Vert ^{2}\leq 2\func{Re}%
\left\langle \left( B^{\ast }A\right) x,x\right\rangle +r^{2}  \label{3.4}
\end{equation}%
for any $x\in H,$ $\left\Vert x\right\Vert =1.$

Since, by (\ref{3.1}),%
\begin{equation*}
\left\Vert Bx\right\Vert ^{2}\geq \frac{1}{\left\Vert B^{-1}\right\Vert ^{2}}%
\left\Vert x\right\Vert ^{2},\quad x\in H
\end{equation*}%
and $\func{Re}\left\langle \left( B^{\ast }A\right) x,x\right\rangle \leq
\left\vert \left\langle \left( B^{\ast }A\right) x,x\right\rangle
\right\vert ,$ hence by (\ref{3.4}) we get%
\begin{equation}
\left\Vert Ax\right\Vert ^{2}+\frac{\left\Vert x\right\Vert ^{2}}{\left\Vert
B^{-1}\right\Vert ^{2}}\leq 2\left\vert \left\langle \left( B^{\ast
}A\right) x,x\right\rangle \right\vert +r^{2}  \label{3.5}
\end{equation}%
for any $x\in H,$ $\left\Vert x\right\Vert =1.$

Taking the supremum over $x\in H,$ $\left\Vert x\right\Vert =1$ in (\ref{3.5}%
), we have%
\begin{equation}
\left\Vert A\right\Vert ^{2}+\frac{1}{\left\Vert B^{-1}\right\Vert ^{2}}\leq
2w\left( B^{\ast }A\right) +r^{2}.  \label{3.6}
\end{equation}%
By the elementary inequality%
\begin{equation}
\frac{2\left\Vert A\right\Vert }{\left\Vert B^{-1}\right\Vert }\leq
\left\Vert A\right\Vert ^{2}+\frac{1}{\left\Vert B^{-1}\right\Vert ^{2}}
\label{3.7}
\end{equation}%
and by (\ref{3.6}) we then deduce the desired result (\ref{3.3}).
\end{proof}

\begin{remark}
If we choose above $B=\lambda I,$ $\lambda \neq 0,$ then we get the
inequality%
\begin{equation}
\left( 0\leq \right) \left\Vert A\right\Vert -w\left( A\right) \leq \frac{1}{%
2\left\vert \lambda \right\vert }r^{2},  \label{3.7.a}
\end{equation}%
provided $\left\Vert A-\lambda I\right\Vert \leq r.$ This result has been
obtained in the earlier paper \cite{SSD2}.

Also, if we assume that $B=\lambda A^{\ast },$ $A$ is invertible, then we
obtain%
\begin{equation}
\left\Vert A\right\Vert \leq \left\Vert A^{-1}\right\Vert \left[ w\left(
A^{2}\right) +\frac{1}{2\left\vert \lambda \right\vert }r^{2}\right] ,
\label{3.7.b}
\end{equation}%
provided $\left\Vert A-\lambda A^{\ast }\right\Vert \leq r,$ $\lambda \neq
0. $
\end{remark}

The following result may be stated as well:

\begin{theorem}
\label{ta.2}Let $A,B:H\rightarrow H$ be two bounded linear operators on $H.$
If $B$ is invertible and for $r>0,$%
\begin{equation}
\left\Vert A-B\right\Vert \leq r,  \label{a.1}
\end{equation}%
then%
\begin{equation}
\left( 0\leq \right) \left\Vert A\right\Vert \left\Vert B\right\Vert
-w\left( B^{\ast }A\right) \leq \frac{1}{2}r^{2}+\frac{\left\Vert
B\right\Vert ^{2}\left\Vert B^{-1}\right\Vert ^{2}-1}{\left\Vert
B^{-1}\right\Vert ^{2}}.  \label{a.2}
\end{equation}
\end{theorem}

\begin{proof}
The condition (\ref{a.1}) is obviously equivalent to%
\begin{equation*}
\left\Vert Ax\right\Vert ^{2}+\left\Vert Bx\right\Vert ^{2}\leq 2\func{Re}%
\left\langle Ax,Bx\right\rangle +r^{2}
\end{equation*}%
for any $x\in H,$ which is clearly equivalent to%
\begin{equation}
\left\Vert Ax\right\Vert ^{2}+\left\Vert B\right\Vert ^{2}\leq 2\func{Re}%
\left\langle B^{\ast }Ax,x\right\rangle +r^{2}+\left\Vert B\right\Vert
^{2}-\left\Vert Bx\right\Vert ^{2}.  \label{a.3}
\end{equation}%
Since%
\begin{equation*}
\func{Re}\left\langle B^{\ast }Ax,x\right\rangle \leq \left\vert
\left\langle B^{\ast }Ax,x\right\rangle \right\vert ,\quad \left\Vert
Bx\right\Vert ^{2}\geq \frac{1}{\left\Vert B^{-1}\right\Vert ^{2}}\left\Vert
x\right\Vert ^{2}
\end{equation*}%
and%
\begin{equation*}
\left\Vert Ax\right\Vert ^{2}+\left\Vert B\right\Vert ^{2}\geq 2\left\Vert
B\right\Vert \left\Vert Ax\right\Vert
\end{equation*}%
for any $x\in H,$ hence by (\ref{a.3}) we get%
\begin{equation}
2\left\Vert B\right\Vert \left\Vert Ax\right\Vert \leq 2\left\vert
\left\langle B^{\ast }Ax,x\right\rangle \right\vert +r^{2}+\frac{\left\Vert
B\right\Vert ^{2}\left\Vert B^{-1}\right\Vert ^{2}-1}{\left\Vert
B^{-1}\right\Vert ^{2}}  \label{a.4}
\end{equation}%
for any $x\in H,$ $\left\Vert x\right\Vert =1.$

Taking the supremum over $x\in H,$ $\left\Vert x\right\Vert =1$ we deduce
the desired result (\ref{a.2}).
\end{proof}

\begin{remark}
If we choose in Theorem \ref{ta.2}, $B=\lambda A^{\ast },$ $\lambda \neq 0,$ 
$A$ is invertible, then we get the inequality:%
\begin{equation}
\left( 0\leq \right) \left\Vert A\right\Vert ^{2}-w\left( A^{2}\right) \leq 
\frac{1}{2\left\vert \lambda \right\vert }r^{2}+\left\vert \lambda
\right\vert \cdot \frac{\left\Vert A\right\Vert ^{2}\left\Vert
A^{-1}\right\Vert ^{2}-1}{\left\Vert A^{-1}\right\Vert ^{2}}  \label{3.9.1}
\end{equation}%
provided $\left\Vert A-\lambda A^{\ast }\right\Vert \leq r.$
\end{remark}

The following result may be stated as well.

\begin{theorem}
\label{tb.1}Let $A,B:H\rightarrow H$ be two bounded linear operators on $H.$
If $B$ is invertible and for $r>0$ we have%
\begin{equation}
\left\Vert A-B\right\Vert \leq r<\left\Vert B\right\Vert ,  \label{b.1}
\end{equation}%
then%
\begin{equation}
\left\Vert A\right\Vert \leq \frac{1}{\sqrt{\left\Vert B\right\Vert
^{2}-r^{2}}}\left( w\left( B^{\ast }A\right) +\frac{\left\Vert B\right\Vert
^{2}\left\Vert B^{-1}\right\Vert ^{2}-1}{2\left\Vert B^{-1}\right\Vert ^{2}}%
\right) .  \label{b.2}
\end{equation}
\end{theorem}

\begin{proof}
The first part of condition (\ref{b.1}) is obviously equivalent to%
\begin{equation*}
\left\Vert Ax\right\Vert ^{2}+\left\Vert Bx\right\Vert ^{2}\leq 2\func{Re}%
\left\langle Ax,Bx\right\rangle +r^{2}
\end{equation*}%
for any $x\in H,$ which is clearly equivalent to%
\begin{equation}
\left\Vert Ax\right\Vert ^{2}+\left\Vert B\right\Vert ^{2}-r^{2}\leq 2\func{%
Re}\left\langle B^{\ast }Ax,x\right\rangle +\left\Vert B\right\Vert
^{2}-\left\Vert Bx\right\Vert ^{2}.  \label{b.3}
\end{equation}%
Since%
\begin{gather*}
\func{Re}\left\langle B^{\ast }Ax,x\right\rangle \leq \left\vert
\left\langle B^{\ast }Ax,x\right\rangle \right\vert , \\
\left\Vert Bx\right\Vert ^{2}\geq \frac{1}{\left\Vert B^{-1}\right\Vert ^{2}}%
\left\Vert x\right\Vert ^{2}
\end{gather*}%
and, by the second part of (\ref{b.1}),%
\begin{equation*}
\left\Vert Ax\right\Vert ^{2}+\left\Vert B\right\Vert ^{2}-r^{2}\geq 2\sqrt{%
\left\Vert B\right\Vert ^{2}-r^{2}}\left\Vert Ax\right\Vert ,
\end{equation*}%
for any $x\in H,$ hence by (\ref{b.3}) we get%
\begin{equation}
2\left\Vert Ax\right\Vert \sqrt{\left\Vert B\right\Vert ^{2}-r^{2}}\leq
2\left\vert \left\langle B^{\ast }Ax,x\right\rangle \right\vert +\frac{%
\left\Vert B\right\Vert ^{2}\left\Vert B^{-1}\right\Vert ^{2}-1}{\left\Vert
B^{-1}\right\Vert ^{2}}  \label{b.4}
\end{equation}%
for any $x\in H,$ $\left\Vert x\right\Vert =1.$

Taking the supremum over $x\in H,$ $\left\Vert x\right\Vert =1$ in (\ref{b.4}%
), we deduce the desired inequality (\ref{b.2}).
\end{proof}

\begin{remark}
The above Theorem \ref{tb.1} has some particular cases of interest. For
instance, if we choose $B=\lambda I,$ with $\left\vert \lambda \right\vert
>r,$ then (\ref{b.1}) is obviously fulfilled and by (\ref{b.2}) we get%
\begin{equation}
\left\Vert A\right\Vert \leq \frac{w\left( A\right) }{\sqrt{1-\left( \frac{r%
}{\left\vert \lambda \right\vert }\right) ^{2}}},  \label{3.15.1}
\end{equation}%
provided $\left\Vert A-\lambda I\right\Vert \leq r.$ This result has been
obtained in the earlier paper \cite{SSD2}.

On the other hand, if in the above we choose $B=\lambda A^{\ast }$ with $%
\left\Vert A\right\Vert \geq \frac{r}{\left\vert \lambda \right\vert }$ \ $%
\left( \lambda \neq 0\right) ,$ then by (\ref{b.2}) we get%
\begin{equation}
\left\Vert A\right\Vert \leq \frac{1}{\sqrt{\left\Vert A\right\Vert
^{2}-\left( \frac{r}{\left\vert \lambda \right\vert }\right) ^{2}}}\left[
w\left( A^{2}\right) +\left\vert \lambda \right\vert \cdot \frac{\left\Vert
A\right\Vert ^{2}\left\Vert A^{-1}\right\Vert ^{2}-1}{2\left\Vert
A^{-1}\right\Vert ^{2}}\right] ,  \label{3.15.2}
\end{equation}%
provided $\left\Vert A-\lambda A^{\ast }\right\Vert \leq r.$
\end{remark}

The following result may be stated as well.

\begin{theorem}
\label{t3.2}Let $A,B$ and $r$ be as in Theorem \ref{t3.1}. Moreover, if%
\begin{equation}
\left\Vert B^{-1}\right\Vert <\frac{1}{r},  \label{3.7.1}
\end{equation}%
then%
\begin{equation}
\left\Vert A\right\Vert \leq \frac{\left\Vert B^{-1}\right\Vert }{\sqrt{%
1-r^{2}\left\Vert B^{-1}\right\Vert ^{2}}}w\left( B^{\ast }A\right) .
\label{3.8}
\end{equation}
\end{theorem}

\begin{proof}
Observe that, by (\ref{3.6}) we have%
\begin{equation}
\left\Vert A\right\Vert ^{2}+\frac{1-r^{2}\left\Vert B^{-1}\right\Vert ^{2}}{%
\left\Vert B^{-1}\right\Vert ^{2}}\leq 2w\left( B^{\ast }A\right) .
\label{3.9}
\end{equation}%
Utlising the elementary inequality%
\begin{equation}
2\frac{\left\Vert A\right\Vert }{\left\Vert B^{-1}\right\Vert }\sqrt{%
1-r^{2}\left\Vert B^{-1}\right\Vert ^{2}}\leq \left\Vert A\right\Vert ^{2}+%
\frac{1-r^{2}\left\Vert B^{-1}\right\Vert ^{2}}{\left\Vert B^{-1}\right\Vert
^{2}},  \label{3.10}
\end{equation}%
which can be stated since (\ref{3.7}) is assumed to be true, hence by (\ref%
{3.9}) and (\ref{3.10}) we deduce the desired result (\ref{3.8}).
\end{proof}

\begin{remark}
If we assume that $B=\lambda A^{\ast }$ with $\lambda \neq 0$ and $A$ an
invertible operator, then, by applying Theorem \ref{t3.2}, we get the
inequality:%
\begin{equation}
\left\Vert A\right\Vert \leq \frac{\left\Vert A^{-1}\right\Vert w\left(
A^{2}\right) }{\sqrt{\left\vert \lambda \right\vert ^{2}-r^{2}\left\Vert
A^{-1}\right\Vert ^{2}}},  \label{3.19.1}
\end{equation}%
provided $\left\Vert A-\lambda A^{\ast }\right\Vert \leq r$ and $\left\Vert
A^{-1}\right\Vert \leq \frac{\left\vert \lambda \right\vert }{r}.$
\end{remark}

The following result may be stated as well.

\begin{theorem}
\label{t3.3}Let $A,B:H\rightarrow H$ be two bounded linear operators. If $%
r>0 $ and $B$ is invertible with the property that $\left\Vert
A-B\right\Vert \leq r$ and%
\begin{equation}
\frac{1}{\sqrt{r^{2}+1}}\leq \left\Vert B^{-1}\right\Vert <\frac{1}{r},
\label{3.11}
\end{equation}%
then%
\begin{equation}
\left\Vert A\right\Vert ^{2}\leq w^{2}\left( B^{\ast }A\right) +2w\left(
B^{\ast }A\right) \cdot \frac{\left\Vert B^{-1}\right\Vert -\sqrt{%
1-r^{2}\left\Vert B^{-1}\right\Vert ^{2}}}{\left\Vert B^{-1}\right\Vert }.
\label{3.12}
\end{equation}
\end{theorem}

\begin{proof}
Let $x\in H,$ $\left\Vert x\right\Vert =1.$ Then by (\ref{3.5}) we have%
\begin{equation}
\left\Vert Ax\right\Vert ^{2}+\frac{1}{\left\Vert B^{-1}\right\Vert ^{2}}%
\leq 2\left\vert \left\langle B^{\ast }Ax,x\right\rangle \right\vert +r^{2},
\label{3.13}
\end{equation}%
and since%
\begin{equation*}
\frac{1}{\left\Vert B^{-1}\right\Vert ^{2}}-r^{2}>0,
\end{equation*}%
we can conclude that $\left\vert \left\langle B^{\ast }Ax,x\right\rangle
\right\vert >0$ for any $x\in H,$ $\left\Vert x\right\Vert =1.$

Dividing in (\ref{3.13}) with $\left\vert \left\langle B^{\ast
}Ax,x\right\rangle \right\vert >0,$ we obtain%
\begin{equation}
\frac{\left\Vert Ax\right\Vert ^{2}}{\left\vert \left\langle B^{\ast
}Ax,x\right\rangle \right\vert }\leq 2+\frac{r^{2}}{\left\vert \left\langle
B^{\ast }Ax,x\right\rangle \right\vert }-\frac{1}{\left\Vert
B^{-1}\right\Vert ^{2}\left\vert \left\langle B^{\ast }Ax,x\right\rangle
\right\vert }.  \label{3.14}
\end{equation}%
Subtracting $\left\vert \left\langle B^{\ast }Ax,x\right\rangle \right\vert $
from both sides of (\ref{3.14}), we get%
\begin{align}
& \frac{\left\Vert Ax\right\Vert ^{2}}{\left\vert \left\langle B^{\ast
}Ax,x\right\rangle \right\vert }-\left\vert \left\langle B^{\ast
}Ax,x\right\rangle \right\vert  \label{3.15} \\
& \leq 2-\left\vert \left\langle B^{\ast }Ax,x\right\rangle \right\vert -%
\frac{1-r^{2}\left\Vert B^{-1}\right\Vert ^{2}}{\left\vert \left\langle
B^{\ast }Ax,x\right\rangle \right\vert \left\Vert B^{-1}\right\Vert ^{2}} 
\notag \\
& =2-\frac{2\sqrt{1-r^{2}\left\Vert B^{-1}\right\Vert ^{2}}}{\left\Vert
B^{-1}\right\Vert }-\left( \sqrt{\left\vert \left\langle B^{\ast
}Ax,x\right\rangle \right\vert }-\frac{\sqrt{1-r^{2}\left\Vert
B^{-1}\right\Vert ^{2}}}{\left\Vert B^{-1}\right\Vert \sqrt{\left\vert
\left\langle B^{\ast }Ax,x\right\rangle \right\vert }}\right) ^{2}  \notag \\
& \leq 2\left( \frac{\left\Vert B^{-1}\right\Vert -\sqrt{1-r^{2}\left\Vert
B^{-1}\right\Vert ^{2}}}{\left\Vert B^{-1}\right\Vert }\right) ,  \notag
\end{align}%
which gives:%
\begin{equation}
\left\Vert Ax\right\Vert ^{2}\leq \left\vert \left\langle B^{\ast
}Ax,x\right\rangle \right\vert ^{2}+2\left\vert \left\langle B^{\ast
}Ax,x\right\rangle \right\vert \frac{\left\Vert B^{-1}\right\Vert -\sqrt{%
1-r^{2}\left\Vert B^{-1}\right\Vert ^{2}}}{\left\Vert B^{-1}\right\Vert }.
\label{3.16}
\end{equation}%
We also remark that, by (\ref{3.11}) the quantity%
\begin{equation*}
\left\Vert B^{-1}\right\Vert -\sqrt{1-r^{2}\left\Vert B^{-1}\right\Vert ^{2}}%
\geq 0,
\end{equation*}%
hence, on taking the supremum in (\ref{3.16}) over $x\in H,$ $\left\Vert
x\right\Vert =1,$ we deduce the desired inequality.
\end{proof}

\begin{remark}
It is interesting to remark that if we assume $\lambda \in \mathbb{C}$ with $%
0<r\leq \left\vert \lambda \right\vert \leq \sqrt{r^{2}+1}$ and $\left\Vert
A-\lambda I\right\Vert \leq r,$ then by (\ref{3.2}) we can state the
following inequality:%
\begin{equation}
\left\Vert A\right\Vert ^{2}\leq \left\vert \lambda \right\vert ^{2}w\left(
A^{2}\right) +2\left\vert \lambda \right\vert \left( 1-\sqrt{\left\vert
\lambda \right\vert ^{2}-r^{2}}\right) w\left( A\right) .  \label{3.25.1}
\end{equation}%
Also, if $\left\Vert A-A^{\ast }\right\Vert \leq r,$ $A$ is invertible and $%
\frac{1}{\sqrt{r^{2}+1}}\leq \left\Vert A^{-1}\right\Vert \leq \frac{1}{r},$
then, by (\ref{3.12}) we also have%
\begin{equation}
\left\Vert A\right\Vert ^{2}\leq w^{2}\left( A^{2}\right) +2w\left(
A^{2}\right) \cdot \frac{\left\Vert A^{-1}\right\Vert -\sqrt{%
1-r^{2}\left\Vert A^{-1}\right\Vert ^{2}}}{\left\Vert A^{-1}\right\Vert }.
\label{3.25.2}
\end{equation}
\end{remark}

One can also prove the following result.

\begin{theorem}
\label{t3.4}Let $A,B:H\rightarrow H$ be two bounded linear operators. If $%
r>0 $ and $B$ is invertible with the property that $\left\Vert
A-B\right\Vert \leq r$ and $\left\Vert B^{-1}\right\Vert \leq \frac{1}{r},$
then%
\begin{align}
(0& \leq )\left\Vert A\right\Vert ^{2}\left\Vert B\right\Vert
^{2}-w^{2}\left( B^{\ast }A\right)  \label{3.17} \\
& \leq 2w\left( B^{\ast }A\right) \cdot \frac{\left\Vert B\right\Vert }{%
\left\Vert B^{-1}\right\Vert }\left( \left\Vert B\right\Vert \left\Vert
B^{-1}\right\Vert -\sqrt{1-r^{2}\left\Vert B^{-1}\right\Vert ^{2}}\right) . 
\notag
\end{align}
\end{theorem}

\begin{proof}
We subtract the quantity $\frac{\left\vert \left\langle B^{\ast
}Ax,x\right\rangle \right\vert }{\left\Vert B\right\Vert ^{2}}$ from both
sides of (\ref{3.14}) to obtain%
\begin{align}
0& \leq \frac{\left\Vert Ax\right\Vert ^{2}}{\left\vert \left\langle B^{\ast
}Ax,x\right\rangle \right\vert }-\frac{\left\vert \left\langle B^{\ast
}Ax,x\right\rangle \right\vert }{\left\Vert B\right\Vert ^{2}}  \label{3.18}
\\
& \leq 2-\frac{\left\vert \left\langle B^{\ast }Ax,x\right\rangle
\right\vert }{\left\Vert B\right\Vert ^{2}}-\frac{1-r^{2}\left\Vert
B^{-1}\right\Vert ^{2}}{\left\vert \left\langle B^{\ast }Ax,x\right\rangle
\right\vert \left\Vert B^{-1}\right\Vert ^{2}}  \notag \\
& =2-2\cdot \frac{\sqrt{1-r^{2}\left\Vert B^{-1}\right\Vert ^{2}}}{%
\left\Vert B\right\Vert \left\Vert B^{-1}\right\Vert }-\left( \frac{\sqrt{%
\left\vert \left\langle B^{\ast }Ax,x\right\rangle \right\vert }}{\left\Vert
B\right\Vert }-\frac{\sqrt{1-r^{2}\left\Vert B^{-1}\right\Vert ^{2}}}{\sqrt{%
\left\vert \left\langle B^{\ast }Ax,x\right\rangle \right\vert }\left\Vert
B^{-1}\right\Vert }\right) ^{2}  \notag \\
& \leq 2\cdot \frac{\left( \left\Vert B\right\Vert \left\Vert
B^{-1}\right\Vert -\sqrt{1-r^{2}\left\Vert B^{-1}\right\Vert ^{2}}\right) }{%
\left\Vert B\right\Vert \left\Vert B^{-1}\right\Vert },  \notag
\end{align}%
which is equivalent with%
\begin{align}
(0& \leq )\left\Vert Ax\right\Vert ^{2}\left\Vert B\right\Vert
^{2}-\left\vert \left\langle B^{\ast }Ax,x\right\rangle \right\vert ^{2}
\label{3.19} \\
& \leq 2\frac{\left\Vert B\right\Vert }{\left\Vert B^{-1}\right\Vert }%
\left\vert \left\langle B^{\ast }Ax,x\right\rangle \right\vert \left(
\left\Vert B\right\Vert \left\Vert B^{-1}\right\Vert -\sqrt{%
1-r^{2}\left\Vert B^{-1}\right\Vert ^{2}}\right)  \notag
\end{align}%
for any $x\in H,$ $\left\Vert x\right\Vert =1.$

The inequality (\ref{3.19}) also shows that $\left\Vert B\right\Vert
\left\Vert B^{-1}\right\Vert \geq \sqrt{1-r^{2}\left\Vert B^{-1}\right\Vert
^{2}}$ and then, by (\ref{3.19}), we get%
\begin{multline}
\left\Vert Ax\right\Vert ^{2}\left\Vert B\right\Vert ^{2}\leq \left\vert
\left\langle B^{\ast }Ax,x\right\rangle \right\vert ^{2}  \label{3.20} \\
+2\frac{\left\Vert B\right\Vert }{\left\Vert B^{-1}\right\Vert }\left\vert
\left\langle B^{\ast }Ax,x\right\rangle \right\vert \left( \left\Vert
B\right\Vert \left\Vert B^{-1}\right\Vert -\sqrt{1-r^{2}\left\Vert
B^{-1}\right\Vert ^{2}}\right)
\end{multline}%
for any $x\in X,$ $\left\Vert x\right\Vert =1.$

Taking the supremum in (\ref{3.20}) we deduce the desired inequality (\ref%
{3.17}).
\end{proof}

\begin{remark}
The above Theorem \ref{t3.4} has some particular instances of interest as
follows. If, for instance, we choose $B=\lambda I$ with $\left\vert \lambda
\right\vert \geq r>0$ and $\left\Vert A-\lambda I\right\Vert \leq r,$ then
by (\ref{3.17}) we obtain the inequality%
\begin{align}
(0& \leq )\left\Vert A\right\Vert ^{2}-w^{2}\left( A\right)  \label{3.29.1}
\\
& \leq 2\left\vert \lambda \right\vert w\left( A\right) \left( 1-\sqrt{1-%
\frac{r^{2}}{\left\vert \lambda \right\vert ^{2}}}\right) .  \notag
\end{align}%
Also, if $A$ is invertible, $\left\Vert A-\lambda A^{\ast }\right\Vert \leq
r $ and $\left\Vert A^{-1}\right\Vert \leq \frac{\left\vert \lambda
\right\vert }{r},$ then by (\ref{3.17}) we can state:%
\begin{align}
(0& \leq )\left\Vert A\right\Vert ^{2}-w^{2}\left( A^{2}\right)
\label{3.29.2} \\
& \leq 2\left\vert \lambda \right\vert w\left( A^{2}\right) \cdot \frac{%
\left\Vert A\right\Vert }{\left\Vert A^{-1}\right\Vert }\left( \left\Vert
A\right\Vert \left\Vert A^{-1}\right\Vert -\sqrt{1-\frac{r^{2}}{\left\vert
\lambda \right\vert ^{2}}\left\Vert A^{-1}\right\Vert ^{2}}\right) .  \notag
\end{align}
\end{remark}

\end{document}